\documentclass[twocolumn]{autart}
                                                          
 \usepackage{cite,comment}
\usepackage{mathrsfs}
\usepackage{graphicx}
\usepackage{graphics}
\usepackage{amssymb}
\usepackage{color}
\usepackage{xspace}
\usepackage{algpseudocode}
\usepackage{bbm}
\usepackage{comment}
\usepackage{algorithmicx}
\usepackage[caption=false]{subfig}
\usepackage{psfrag}
\usepackage{url}
\usepackage{float}
\usepackage{caption}
\usepackage{soul}

\usepackage{rotating}
\usepackage{chngcntr}
\usepackage{apptools}
\AtAppendix{\counterwithin{thm}{section}}
\usepackage{graphicx}
\usepackage{graphics}
\usepackage{amssymb}
\usepackage{amsmath}
\usepackage{mathtools}
\usepackage{color}
\usepackage{xspace}
\usepackage{algpseudocode}
\usepackage{bbm}
\usepackage{comment}
\usepackage{algorithmicx}
\usepackage{subfig}
\usepackage{psfrag}
\usepackage{tikz}
\usetikzlibrary{shapes,arrows}
\usetikzlibrary{positioning}

\usepackage{rotating}
\usepackage{chngcntr}
\usepackage{apptools}
\usepackage{listings}
\AtAppendix{\counterwithin{thm}{section}}
\usepackage{graphicx}
\usepackage{graphics}
\usepackage{amssymb}
\usepackage{amsmath}
\usepackage{color}
\usepackage{xspace}
\usepackage{algpseudocode}
\usepackage{bbm}
\usepackage{comment}
\usepackage{algorithmicx}
\usepackage{subfig}
\usepackage{psfrag}
\usepackage{tikz}
\usetikzlibrary{shapes,arrows}
\usetikzlibrary{positioning}

\usepackage{enumitem}
\usepackage{graphicx}               
\usepackage{amsmath,amssymb}

\usepackage{graphics}
\usepackage{amssymb}
\usepackage{amsmath}
\usepackage{color}
\usepackage{xspace}
\usepackage{algpseudocode}
\usepackage{bbm}
\usepackage{comment}
\usepackage{algorithmicx}
\usepackage{subfig}
\usepackage{psfrag}
\usepackage{soul}
\usepackage{tikz}
\usetikzlibrary{petri}
\usetikzlibrary{shapes,arrows}
\usetikzlibrary{positioning}

\usepackage{cancel}

\tikzstyle{block}=[draw opacity=0.7,line width=1.4cm]

\usepackage{accents}
\newcommand{\ubar}[1]{\underaccent{\bar}{#1}}
	
\newcommand{\oprocendsymbol}{\hbox{$\bullet$}}
\newcommand{\oprocend}{\relax\ifmmode\else\unskip\hfill\fi\oprocendsymbol}

\newcommand{\VV}{\mathcal{V}}
\newcommand{\EE}{\mathcal{E}}
\newcommand{\GG}{\mathcal{G}}

\newcommand{\real}{{\mathbb{R}}}
\newcommand{\reals}{{\mathbb{R}}}
\newcommand{\realpositive}{{\mathbb{R}}_{>0}}

\newcommand{\rank}{\operatorname{rank}}

\newcommand{\until}[1]{\in\{1,\dots,#1\}}
\newcommand{\vect}[1]{\boldsymbol{\mathbf{#1}}}
\newcommand{\vectsf}[1]{\boldsymbol{\mathbf{\mathsf{#1}}}}

\newcommand{\Tvect}[1]{\tilde{\boldsymbol{\mathbf{#1}}}}

\newcommand{\Hvect}[1]{\hat{\boldsymbol{\mathbf{#1}}}}

\newcommand{\dvect}[1]{\dot{\vect{#1}}}

\newcommand{\ee}{\operatorname{e}}

\newcommand{\Diag}[1]{\operatorname{Diag}(#1)}

 \newcommand{\boxend}{\hfill \ensuremath{\Box}}

\newtheorem{thm}{Theorem}[section]

\newtheorem{rem}{Remark}[section]

\newtheorem{lem}{Lemma}[section]

\newtheorem{assump}{Assumption}[section]


\tikzstyle{block}=[draw opacity=0.7,line width=1.4cm]
\DeclareMathAlphabet{\mathpzc}{OT1}{pzc}{m}{it}
\definecolor{CranJ}{cmyk}{0,0.69,0.54,0.04} 
\definecolor{PinkJ}{cmyk}{0,0.71,0.43,0.12} 
\definecolor{Cran}{cmyk}{0,0.73,0.41,0.29} 
\definecolor{VRed}{cmyk}{0,0.75,0.25,0.2} 
\definecolor{ORed}{cmyk}{0,0.75,0.75,0} 
\definecolor{CBlue}{cmyk}{1,0.25,0,0} 

\begin{document}
\begin{frontmatter}
   \runtitle{Distributed optimal resource allocation}
  
   \title{{\large A Distributed Continuous-time Modified Newton-Raphson Algorithm }
   } 
  
  \thanks[footnoteinfo]{Corresponding author: H. Moradian. This work is supported by NSF award ECCS-1653838.}
      \author[Paestum]{Hossein Moradian}\ead{hmoradia@uci.edu}\quad
  \author[Paestum]{Solmaz S. Kia}\ead{solmaz@uci.edu}

  \address[Paestum]{Department of Mechanical and Aerospace
    Engineering, University of California, Irvine}
  \begin{abstract}
We propose a continuous-time second-order optimization algorithm for solving unconstrained convex optimization problems with bounded Hessian. We show that this alternative algorithm has a comparable convergence rate to that of the continuous-time Newton-Raphson method, however structurally, it is amenable to a more efficient distributed implementation. We present a distributed implementation of our proposed optimization algorithm and prove its convergence via Lyapunov analysis. A numerical example demonstrates our results. 
  \end{abstract}
\end{frontmatter}

\pagenumbering{roman}
\maketitle

\section{INTRODUCTION}
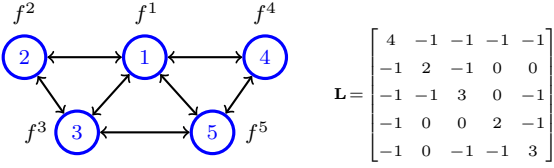
\begin{figure}[!t]
  \unitlength=0.5in \centering
\begin{center}
   \begin{tikzpicture}[auto,thick,scale=1, every node/.style={scale=1}]
     \node (1) at (0,1) [ draw, minimum size=7pt,color=blue, circle, very thick, label=above:{\scriptsize$f^1\!$}] {{\scriptsize 1}};
      \node (2) at (-1.6,1) [ draw, minimum size=7pt,color=blue, circle,very thick, label=above:{\scriptsize$f^2$}] {{\scriptsize 2}};
       \node (3) at (-0.9,0) [ draw, minimum size=7pt,color=blue, circle, very thick, label=left:{\scriptsize$f^3\!$}] {{\scriptsize 3}};
    \node (4) at (1.6,1) [ draw, minimum size=7pt,color=blue, circle, very thick, label=above:{\scriptsize$f^4$}] {{\scriptsize 4}};
    \node (5) at (0.9,0) [ draw, minimum size=7pt,color=blue, circle, very thick, label=right:{\scriptsize$f^5$}] {{\scriptsize 5}};
    \node (L) at (4,0.5) [  minimum size=7pt,color=black] {{\tiny $\vect{L}\!=\!\begin{bmatrix}4&-1&-1&-1&-1\\-1&2&-1&0&0\\-1&-1&3&0&-1\\-1&0&0&2&-1\\-1&0&-1&-1&3\end{bmatrix}$}};
     \path[draw, <->]  (1)--(2) ;
           \path[draw, <->] (1)--(3);
           \path[draw, <->] (2)--(3);
      \path[draw, <->] (3)--(5) ;
     \path[draw, <->] (1)--(4) ;
     \path[draw, <->] (4)--(5) ;
     \path[draw, <->]  (5)--(1) ;
\end{tikzpicture}
\end{center}
        \caption{{ A \emph{connected graph}: in a connected undirected graph agents connected by an edge can exchange information. Moreover, there is a path from any agent to any other agent. }}
    \label{fig:network}
\end{figure}

Consider a network of $N$ agents interacting over a  connected graph~$\GG$, see Fig.~\ref{fig:network}. Each agent
$i \until{N}$ is endowed with a local cost function
$f^i:\mathbb{R}^d\to\mathbb{R}$ which is  twice differentiable and $m^i$-strongly convex. Our objective is to
design a distributed optimization algorithm such that each agent
obtains the global minimizer $\vect{x}^\star\in\real^d$
of the feasible optimization problem
\vspace{-0.1in}
\begin{equation}\label{eq::glob_opt_prob}
  \vect{x}^{\star} = \arg\underset{{\vect{x}\in
      \reals^d}}{\min} 
  \,\,f(\vect{x}), \quad f(\vect{x}) = \sum\nolimits_{i=1}^N f^i(\vect{{x}}),\vspace{-0.1in}
\end{equation}
using local interactions with its neighbors. The existing distributed optimization solutions are mostly consensus-based approaches that use gradient and sub-gradient methods, see e.g.,~\cite{AN-AO:09,BJ-MR-MJ:09,SB-NP-EC-BP-JE:10,MZ-SM:09c,JD-AA-MW:12} for discrete-time and~\cite{JW-NE:11,FZ-DV-AC-GP-LS:11,JL-CYT:12,SSK-JC-SM:15-auto,GD-HK-ME:14}
for continuous-time algorithms. Even though the  gradient-based solutions' distributed implementation is fully understood and requires low computational resources, they suffer from a low convergence rate, especially near the solution. With the recent advances in fast computing via graphics processing units (GPUs), the interest in Newton-based optimization algorithms, which use second-order information to achieve faster convergence, is renewed for large-scale optimization problems~\cite{ZY-AG-SS-MM-KK-MWM:20,JFH-SE-SA-AV:18}. The popular Newton-Raphson (NR) method uses the inverse of the Hessian of the total cost multiplied by the gradient of the total cost, i.e., $-(\sum\nolimits_{i=1}^N\vect{H}^i(\vect{x}))^{-1}(\sum\nolimits_{i=1}^N\vect{g}^i(\vect{x}))$, as the descent direction. Here, $\vect{g}^i(\vect{x})=\nabla f^i(\vect{x})$ and $\vect{H}^i(\vect{x})=\nabla^2 f^i(\vect{x})$. Starting from a local guess $\vect{x}^i\in\real^d$, $i\until{N}$, a common way to execute the NR algorithm in a decentralized way is to use a consensus-based framework to track  $\sum\nolimits_{i=1}^N\vect{H}^i(\vect{x}^i)$ and $\sum\nolimits_{i=1}^N\vect{g}^i(\vect{x}^i)$ for every agent cooperatively
 by exchanging the gradient and the Hessian of the local costs, see e.g.,~\cite{DV-FZ-AC-GP-LS:15} for continuous-time and~\cite{DV-FZ-AC-GP-LS:15,NB-RC-GN-LS-DV:19}  for discrete-time algorithms. These algorithms result in $O(Nd^2)$
  communication, computation and storage costs per agent to solve problem~\eqref{eq::glob_opt_prob}. To remove communicating Hessian among agents,~\cite{AM-QL-AR:15} proposes a distributed algorithm that approximates Newton step by truncating the Taylor series expansion of the exact Newton step at $K$ terms. But, implementing this algorithm requires aggregating information from $K$ hops away. Increasing $K$ makes the method arbitrarily close to Newton’s method at the cost of increasing the communication overhead of each iteration. Building on~\cite{AM-QL-AR:15},~\cite{FM-EW:20} proposes an asynchronous implantation to manage communication cost but the method works only for univariate local cost functions.

 In this paper, we provide an alternative continuous-time second-order algorithm with a comparable convergence rate to that of the continuous-time NR algorithm, but with a structure that is amenable to a more resource-efficient distributed implementation that also requires information exchange only among one-hop neighbors. In the distributed implementation of this algorithm, agents use the inverse of their local Hessians but do not need to communicate it. As a result, our proposed algorithm's communication, computation, and storage cost per agent are $O(Nd)$.
 We establish the convergence of our algorithm using Lyapunov stability analysis. Simulations demonstrate our results.

\section{Preliminaries}\label{sec::prelim}
Our notations are standard and definitions are given if it is necessary to avoid confusion. A differentiable function $f: \reals^d\to\reals$ is \emph{$m$-strongly convex} ($m\in\real_{>0}$) in a set $C$ if and only if $(\vect{z}-\vect{x})^\top(
  \nabla f(\vect{z})-\nabla f(\vect{x}))\geq
  m\|\vect{z}-\vect{x}\|^2, ~\forall
  ~\vect{x},\vect{z}\in C,~\vect{x}\neq\vect{z}$.
For twice differentiable function $f$ the $m$-strong convexity  ($m>0$) is also equivalent to $\vect{H}(\vect{x})=\nabla^2 f(\vect{x})\geq
  m \vect{I}$,  $\forall
  ~\vect{x}\in C.$

 A connected graph, see Fig.~\ref{fig:network}, is represented by $\GG =
(\VV ,\EE ,\vect{{A}})$, where $\VV=\{1,\dots,N\}$ is the \emph{node set},
$\EE=\{e_1,\cdots,e_M\} \subseteq \VV\times \VV$ is the \emph{edge set}, and $\vect{{A}}=[a_{ij}]\in\real^{N\times N}$ is the \emph{adjacency} matrix such that $a_{ij}=a_{ji}=1$ if $(i, j) \in\EE$ and $
{a}_{ij} = 0$, otherwise. The \emph{incidence} matrix is  $\vect{B}=\frac{1}{\sqrt{2}}[\vect{b}^1,\,-\vect{b}^1,\cdots,\vect{b}^M,\,-\vect{b}^M]\in\real^{N\times 2M}$  where $\vect{b}^k\in\real^N$ is a vector corresponds to the edge $e_k=(i,j)\in\EE$ with zero elements except for $i$th and $j$th components with respectively $b_i=1$, $b_j=-1$. 
The \emph{Laplacian} matrix of a graph is $\vect{L} =
\Diag{\vect{{A}} \vect{1}_N}- \vect{{A}}$. Note that $\vect{L}\vect{1}_N=\vect{0}$. Moreover, $\vect{L}=\vect{B}\vect{B}^\top$.  A graph is connected
 if and only if $\vect{1}_N^T\vect{L}=\vect{0}$,  and $\rank(\vect{L})=N-1$. 
For a connected graph, eigenvalue of $\vect{L}$ are $\lambda_1=0$, $\{\lambda_i\}_{i=1}^N\subset\real_{>0}$. We let  $\lambda_i\leq \lambda_j$, for $i<j$. Moreover, 
\begin{align}\label{eq::L_inverse}
\vect{\Pi}_N=\vect{L}\vect{L}^+=\vect{B}\vect{B}^\top(\vect{B}\vect{B}^\top)^+=\vect{B}(\vect{B}^\top\vect{B})^+\vect{B}^\top\!,
\end{align}
 
where $(.)^+$ denotes the generalized inverse matrix~\cite{JG-RAF:19} and $\vect{\Pi}_N= \vect{I} - \frac{1}{N}\vect{1}_N\vect{1}_N^\top$, where $\vect{1}_N$ is the vector of $N$~ones.

Throughout the paper, the following assumption holds.
\begin{assump}\label{asm:convexity}\rm{
The local cost functions $f^i:\real^d\rightarrow\real$ are $m^i$-strongly convex with the bounded Hessians $\ubar{m}^i\vect{I}\leq\vect{H}^i(\vect{x}^i)\leq\bar{m}^i\vect{I}$, for some $\ubar{m}^i,\bar{m}^i\in\real_{>0}$.}\boxend
\end{assump}
Thus, the total cost $f$ is $\ubar{m}$-strongly convex and its~Hessian is upper-bounded by $\bar{m}\vect{I}$. Moreover, $\vect{x}^\star$ in~\eqref{eq::glob_opt_prob} is unique~\cite{DPB:99}. Here, 
$\ubar{m}\!=\!\min\{\ubar{m}^i\}_{i=1}^N,~ \bar{m}\!=\!\max\{\bar{m}^i\}_{i=1}^N$.

\section{Problem Definition}\label{sec::Prob_formu}
The continuous-time NR algorithm to solve problem~\eqref{eq::glob_opt_prob} is given by
\begin{align}\label{eq::CHessSolv}
\text{NR:}~~~~~&\dvect{x}=-\Big(\sum\nolimits_{i=1}^N\vect{H}^i(\vect{x})\Big)^{-1}\sum\nolimits_{i=1}^N\vect{g}^i(\vect{x}).
\end{align}
In this paper, we propose the Hessian inverse sum optimization (HISO) algorithm
\begin{equation}\label{eq::CHessSolv2} 
\text{HISO:}~~~~\dvect{x}=-\Big(\frac{1}{N}\sum\nolimits_{i=1}^N\vect{H}^i(\vect{x})^{-1}\Big)\sum\nolimits_{i=1}^N\vect{g}^i(\vect{x}),
\end{equation} 
as an alternative Hessian-based solver for the optimization problem~\eqref{eq::glob_opt_prob}. We show that this algorithm has the convergence rate no worse than that of~\eqref{eq::CHessSolv} but it has a structure that is amenable to a distributed implementation with more efficient resource usage. We start by the auxiliary result below.

\begin{lem}(Bound on the inverse of sum of  symmetric positive definite matrices)\label{lem::PD_product} {\rm Let every $\vect{H}^i\in\real^{d\times d}$, $i\in\{1,\cdots,N\}$, be a positive definite matrix. Then }
\begin{align}\label{eq::inequ_hessian}
\Big(\sum\nolimits_{i=1}^N\vect{H}^i\Big)^{-1}\leq\frac{1}{N}\sum\nolimits_{i=1}^N{\vect{H}^i}^{-1}.
\end{align}
\end{lem}
\begin{pf}
The proof is by mathematical induction. 
Recall that the inverse of positive definite matrices is a convex function~\cite{KN:11}. Hence, for $N=2$ for any $\kappa\in[0,1]$ we have
\begin{align}\label{eq::cov_inverse}
    \big(\kappa\,\vect{H}^1+(1-\kappa)\,\vect{H}^2\big)^{-1}\leq\kappa\,{\vect{H}^1}^{-1}+(1-\kappa)\,{\vect{H}^2}^{-1}.
\end{align}
Substituting $\kappa=0.5$ gives $(\vect{H}^1+\vect{H}^2)^{-1}\leq\frac{1}{4}({\vect{H}^1}^{-1}+{\vect{H}^2}^{-1})\leq\frac{1}{2}({\vect{H}^1}^{-1}+{\vect{H}^2}^{-1})$. Thus,~\eqref{eq::inequ_hessian} holds for $N=2$. 
Next, assuming $(\sum_{i=1}^{N-1}\vect{H}^i)^{-1}\leq\frac{1}{N-1}\sum_{i=1}^{N-1}{\vect{H}^i}^{-1}$ we show that $(\sum_{i=1}^{N}\vect{H}^i)^{-1}\leq\frac{1}{N}\sum_{i=1}^{N}{\vect{H}^i}^{-1}$ holds. To this aim, notice that given~\eqref{eq::cov_inverse} we obtain
\begin{align*}
    &\Big(\frac{1}{N}\vect{H}^N\!\!+\!\frac{N-1}{N}\sum_{i=1}^{N-1}\!\vect{H}^i\Big)^{-1}\!\!\!\leq\!\frac{1}{N}{\vect{H}^N}^{-1}\!\!\!+\!\frac{N-1}{N}\sum_{i=1}^{N-1}\!{\vect{H}^i}^{-1}\\
    &\leq \frac{1}{N}{\vect{H}^N}^{-1} \!\!+\frac{N-1}{N}\Big(\frac{1}{N-1}\sum_{i=1}^{N-1}{\vect{H}^i}^{-1}\Big)=\frac{1}{N}\sum_{i=1}^{N}\!{\vect{H}^i}^{-1}.
\end{align*}
Since $\big(\sum_{i=1}^N\vect{H}^i\big)^{-1}\leq\big(\frac{1}{N}\vect{H}^N+\frac{N-1}{N}\sum_{i=1}^{N-1}(\vect{H}^i)\big)^{-1}$, then~\eqref{eq::inequ_hessian} holds for any $N\geq2$, which concludes proof.
\end{pf}

 Lemma~\ref{lem::PD_product} enables us to make the following statement about the HISO algorithm's convergence guarantees.  
\begin{thm}(Convergence analysis of the HISO algorithm)\label{lem::central_solution}
{\rm
Consider the optimization problem~\eqref{eq::glob_opt_prob} and let Assumption~\ref{asm:convexity} hold.  Then, starting from any initial condition $\vect{x}(0)\in\real^d$, as $t\to\infty$ the HISO algorithm~\eqref{eq::CHessSolv2} converges exponentially fast to $\vect{x}^\star\in\real^d$, the unique minimizer of the optimization problem~\eqref{eq::glob_opt_prob}. Furthermore, the rate of convergence of~\eqref{eq::CHessSolv2} is no worse than the rate of convergence of the algorithm~\eqref{eq::CHessSolv}.}
\end{thm}
\begin{pf}
Consider the candidate Lyapunov function  $V(\vect{x})=f(\vect{x})-f(\vect{x}^\star)$. Given Assumption~\ref{asm:convexity}, we have $\ubar{m}\|\vect{x}-\vect{x}^\star\|^2\leq V(\vect{x})\leq\bar{m}\|\vect{x}-\vect{x}^\star\|^2$, and $\|(\sum\nolimits_{i=1}^N\vect{H}^i(\vect{x}))^{-1}\nabla f(\vect{x})\|\leq\frac{\bar{m}}{\ubar{m}}\|\vect{x}-\vect{x}^\star\|$ where $\nabla f(\vect{x})=\sum\nolimits_{i=1}^N\vect{g}^i(\vect{x})$.
The derivative of
$V(\vect{x})$ along trajectories of~\eqref{eq::CHessSolv2},
satisfies $\dot{V}(x)\!=\!-\frac{1}{N}\nabla f(\vect{x})^\top\!(\sum\nolimits_{i=1}^N\!\!\vect{H}^i(\vect{x})^{-1})\nabla f(\vect{x})\\\leq-\frac{\bar{m}^2}{\ubar{m}}\|\vect{x}-\vect{x}^\star\|^2$, which by virtue of~\cite[Theorem 4.10]{HKK:02} confirms the exponential stability of~\eqref{eq::CHessSolv2}. On the other hand, derivative of $V(\vect{x})$ along~\eqref{eq::CHessSolv} satisfies
$
\dot{V}=-\nabla f(\vect{x})^\top (\sum_{i=1}^N\vect{H}^i(\vect{x}))^{-1}\nabla f(\vect{x})\leq-\frac{\bar{m}^2}{\ubar{m}}\|\vect{x}-\vect{x}^\star\|^2$, confirming the exponential stability of~\eqref{eq::CHessSolv}.
By virtue of Lemma~\ref{lem::PD_product}, 
$\label{eq::proof_int}
    -\frac{1}{N}\nabla f(\vect{x})^\top (\sum\nolimits_{i=1}^N\vect{H}^i(\vect{x})^{-1})\nabla f(\vect{x})\leq-\nabla f(\vect{x})^\top(\sum\nolimits_{i=1}^N\vect{H}^i(\vect{x}))^{-1}\nabla f(\vect{x}),$
   which indicates that the derivative of $V(\vect{x})$ is more negative along the trajectories of~\eqref{eq::CHessSolv2} than~those of~\eqref{eq::CHessSolv}. Hence,
we can conclude that the rate of convergence of~\eqref{eq::CHessSolv2} is no worse than the convergence rate of~\eqref{eq::CHessSolv}.
\end{pf}

\begin{figure}[t]
  \centering
     \includegraphics[scale=0.58]{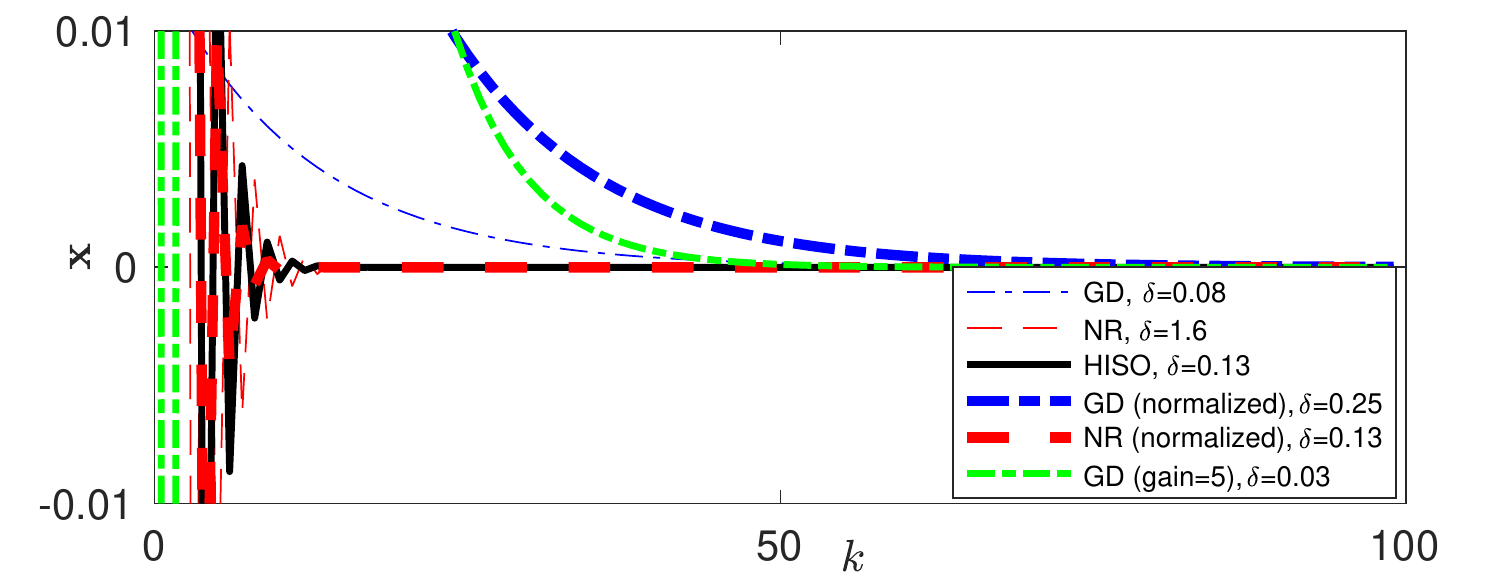}
      \caption{{The  convergence of the Euler discretized GD, NR and HISO algorithms under different conditions to find the minimum of the cost function, $f(\vect{x})=\sum_{i=1}^{10} a^ix^2+b^ix^4$ where $a^i$ and $b^i$ are randomly chosen in $[0,0.1]$. For each algorithm the stepsize is set to its optimum value, obtained numerically,  corresponding to its fastest convergence. 
      }}\label{fig:ex2}
\end{figure}

Comparing rate of convergence of continuous-time
optimization algorithms is a rather subtle matter. Any claim for an accelerated convergence by an algorithm meets the counter-argument that the `simple' continuous-time gradient descent algorithm can be made arbitrarily fast using large scalar multiplicative gains. To address this dilemma, 
one can think of continuous-time algorithms as first-order integrator dynamics $\dot{\vect{x}}=\alpha\,\vect{u}$ with   $\alpha\in\real_{>0}$, where the system input $\alpha\,\vect{u}$ is the control effort of the algorithm. Suppose the control effort is bounded as $\|\alpha\vect{u}\|\leq \alpha\kappa_0\|\vect{x}-\vect{x}^\star\|$,  with $\kappa_0\in\realpositive$. For an exponentially convergent algorithm, by virtue of~\cite[Theorem 4.14]{HKK:02}, there exists a Lyapunov function that satisfies $\kappa_1\|\vect{x}-\vect{x}^\star\|^2\leq V(\vect{x})\leq\kappa_2\|\vect{x}-\vect{x}^\star\|^2$ , and  $\dot{V}\leq-\alpha\kappa_3\|\vect{x}-\vect{x}^\star\|^2$ for some  $\kappa_1,\kappa_2,\kappa_3\in\realpositive$. 
Then, by virtue of~\cite[Theorem 4.10]{HKK:02} the exponential rate of convergence of the algorithm is $\frac{\alpha\kappa_3}{\kappa_2}$, indicating that increasing $\alpha$ increases the rate of convergence. Now, on the other hand, for the Euler-discretized form of the algorithm, i.e., $\vect{x}(k+1)=\vect{x}(k)+\delta\alpha\vect{u}$, $k\in\mathbb{Z}_{\geq0}$, 
using the same Lyapunov function, we obtain $\Delta V(\vect{x}(k))=V(\vect{x}(k+1))-V(\vect{x}(k))\leq-\delta \alpha\kappa_3\|\vect{x}-\vect{x}^\star\|^2+\frac{\delta^2\alpha^2}{2}\vect{u}^\top\nabla^2V(\zeta)\vect{u}$ where, $\zeta\in[\vect{x}(k),\vect{x}(k+1))$. Let  $ \nabla^2V(\zeta)\leq \beta \vect{I}$, which is normally satisfied in optimization problems. Then, we can write $\Delta V(\vect{x}(k))\leq-\delta \alpha(\kappa_3-\frac{\delta \alpha}{2}\beta\kappa_0^2)\|\vect{x}-\vect{x}^\star\|^2$, which indicates that an admissible stepsize $\delta$ for Euler-discretized form of the algorithm should satisfy $0<\delta <\frac{2\kappa_3}{\alpha \beta\kappa_0^2}$. 
Thus, increasing $\alpha$ results in smaller stepsizes. Moreover, algorithms that employ larger control effort (larger $\alpha\,\kappa_0$) will have smaller stepsize.
As such, to be mindful of practical Euler-discretize implementation of continuous-time algorithms, any claim 
to a continuous-time algorithm being faster than another 
should be evaluated under the requirement that the algorithms employ the same maximum control effort level.  
Figure~\ref{fig:ex2} shows the convergence behavior of the discrete-time implementation of the gradient descent (GD), NR, and HISO algorithms. As we can see, NR and HISO show comparable responses and also faster convergence than GD. For all three cases, the maximum control effort happens at the initial time, with GD having the largest and NR having the smallest values. If we normalize the control efforts of NR and GD with respect to that of HISO by using, respectively, gains  $\|\frac{1}{N}\sum\nolimits_{i=1}^N\vect{H}^i(\vect{x}(0))^{-1}\|$ and  $\|\sum\nolimits_{i=1}^N\vect{H}^i(\vect{x}(0))\|\|\frac{1}{N}\sum\nolimits_{i=1}^N\vect{H}^i(\vect{x}(0))^{-1}\|$, the GD algorithm can use larger stepsize and NR should use a smaller stepsize, with GD still showing slower convergence. Notice that, as Fig.~\ref{fig:ex2} shows, if we increase the GD algorithm's control effort by using the gain $\alpha=5$, the discrete-time implementation still has slower convergence because we are forced to use a smaller Euler discretization stepsize.

HISO algorithm uses the sum of the inverse of the Hessian of the local cost functions rather than the inverse of the sum of the local Hessians as in the NR algorithm. This trait, as shown below, results in a more efficient distributed implementation for algorithm~\eqref{eq::CHessSolv2}, in which agents only incur a cost of $O(Nd)$ in communication, computations, and storage rather than $O(Nd^2)$ as in the distributed NR algorithms in the literature~\cite{DV-FZ-AC-GP-LS:15,NB-RC-GN-LS-DV:19}. 

\section{Distributed HISO Algorithm}
\label{sec::main} 
Our proposed distributed implementation of the HISO algorithm is 
\begin{subequations}\label{eq::dis_solver}
\begin{align}
    \vect{z}^i&=\vect{g}^i(\vect{x}^i)+\vect{v}^i,\label{eq::dis_solver_a}\\
    \dvect{v}^i&\!=\!-\!\sum\nolimits_{j=1}^N\!\!a_{ij}\,\text{sgn}(\vect{z}^i\!-\!\vect{z}^j)+\!\sum\nolimits_{j=1}^N\!\!a_{ij}(\vect{x}^i\!-\!\vect{x}^j),\label{eq::dis_solver_b}\\
    \dvect{x}^i&\!=\!-\vect{H}^i(\vect{x})^{-1}\big(\vect{z}^i\!+\! \sum\nolimits_{j=1}^N\!\!a_{ij}(\vect{x}^i-\vect{x}^j)\big)
    ,\label{eq::dis_solver_c}
\end{align}
\end{subequations}
$i\!\until{N}$. Conceptually, our approach to construct~\eqref{eq::dis_solver} was to use the finite-time dynamic average consensus algorithm of~\cite{FC-YC-WR:12} (\eqref{eq::dis_solver_a} and \eqref{eq::dis_solver_b}) with input $\vect{g}^i(\vect{x}^i)\!=\!\nabla f^i(\vect{x}^i)$ to generate $\vect{z}^i\!\rightarrow\!\frac{1}{N}\sum_{j=1}^{N}\!\vect{g}^j(\vect{x}^j)$ as $t\!\to\!\infty$. For algorithm of~\cite{FC-YC-WR:12} to converge we need $\sum_{i=1}^N\!\vect{v}^i(0)=\vect{0}$, which can trivially be satisfied using $\vect{v}^i(0)=\vect{0}$. 
Next, we noticed that the collective dynamics exhibits  $\sum_{i=1}^N\!\dvect{g}^i(\vect{x}^i)\!\to-\sum\nolimits_{j=1}^N\!\vect{g}^j(\vect{x}^j)$ as $\vect{z}^i$ converges. Then, if agreement occurs, every agent has a copy of the HISO algorithm locally. We added $\sum\nolimits_{j=1}^N\!a_{ij}(\vect{x}^i\!-\!\vect{x}^j)$ to~\eqref{eq::dis_solver_b} and~\eqref{eq::dis_solver_c} for technical reasons to create agreement between the decision vector of the~agents. 

In what follows, we provide a formal proof of convergence and stability analysis of~\eqref{eq::dis_solver}.  For analysis, we write algorithm~\eqref{eq::dis_solver} in the compact form
\begin{subequations}\label{eq::dis_solver_com}
\begin{align}
    {\dvect{z}}&=- \,\vectsf{B}\,\text{sgn}(\vectsf{B}^\top\vect{z})+\,\vectsf{B}\vectsf{B}^\top\vect{x}+\frac{\text{d}}{\text{d}t}\vect{g}(\vect{x}),\label{eq::dis_solver_c_b}\\
      \dvect{x}&=-\vect{\mathcal{H}}^{-1}\big(\vect{z}+\,\vectsf{B}\vectsf{B}^\top\vect{x}\big)
      ,\label{eq::dis_solver_c_d}
\end{align}
\end{subequations}
where $\vect{\mathsf{B}}=\vect{B}\otimes\vect{I}_d$, $\vect{\mathcal{H}}= \text{diag}(\vect{H}^1(\vect{x}),\cdots,\vect{H}^{N}(\vect{x}))$, $\vect{g}(\vect{x})=[\vect{g}^1(\vect{x}^1)^\top,\cdots,\vect{g}^N(\vect{x}^N)^\top]^\top$
and $\vect{\mathsf{\Pi}}=\vect{\Pi}_N\otimes\vect{I}_d$ with the network aggregated variables $\vect{z},\vect{x}\in\real^{dN}$. 
The following result shows that agents  arrive at agreement in their $\{\vect{z}^i\}_{i=1}^N$ in finite time, and in their $\{\vect{x}^i\}_{i=1}^N$ as $t\rightarrow\infty$.

\begin{lem}(Consensus  in algorithm~\eqref{eq::dis_solver} over connected graphs) \label{thm::convergence}
{\rm Let $\GG$ be a connected graph. Under  Assumption~\ref{asm:convexity}, starting algorithm~\eqref{eq::dis_solver} over $\GG$ from any $\vect{x}^i(0),\vect{v}^i(0)\in\real^d$, $\sum_{i=1}^N\vect{v}^i(0)=\vect{0}$, every $\vect{x}^i(t)$, $i\until{N}$, converges  to $\frac{1}{N}\sum_{j=1}^N\vect{x}^j$ as time goes to infinity, while every $\vect{z}^i$, $i\until{N}$, converges to $\frac{1}{N}\sum_{j=1}^N\vect{g}^j(\vect{x}^j)$ in finite~time.
}
\end{lem}
\begin{pf}
~\eqref{eq::dis_solver_b} leads to $\sum\nolimits_{i=1}^N\dvect{v}^i=\vect{0}$, which along with  $\sum_{i=1}^N\vect{v}^i(0)=\vect{0}$ gives $\sum_{i=1}^N\vect{v}^i(t)=\vect{0}$ for any $t\in\real_{\geq0}$. Moreover, from~\eqref{eq::dis_solver_a}, we obtain \begin{equation}\label{eq::z_g_sum}\sum\nolimits_{i=1}^N\vect{z}^i(t)=\sum\nolimits_{i=1}^N\vect{g}^i(\vect{x}^i(t)),~~~t\in\real_{\geq0}.\end{equation}
Next, note that $\frac{\text{d}}{\text{d}t}\vect{g}(\vect{x})=\vect{\mathcal{H}}\,\dvect{x}$. Then, we can obtain from~\eqref{eq::dis_solver_com} that $
    \dvect{z}=- \,\vectsf{B}\,\text{sgn}(\vectsf{B}^\top\vect{z})-\vect{z}$,
which gives
\begin{equation}\label{eq::z_sum}\sum_{i=1}^N\dvect{z}^i=-\sum_{i=1}^N\vect{z}^i\rightarrow\sum_{i=1}^N\vect{z}^i(t)=\text{e}^{-t}\sum_{i=1}^N\vect{z}^i(0) , ~~t\in\real_{\geq0}.\end{equation}
Moreover, using 
\begin{subequations}\label{eq::change_of_var}
\begin{align}
 \tilde{\vect{z}}(t)&=\vectsf{\Pi}\,\vect{z}(t),\label{eq::change_of_var_a}\\
 \tilde{\vect{x}}(t)&=\vectsf{\Pi}\,\vect{x}(t),\label{eq::change_of_var_b}
\end{align}
\end{subequations}
~\eqref{eq::dis_solver_com} can be written as
\begin{subequations}\label{eq::new_var}
\begin{align}
   \dot{\Tvect{z}}&\!=- \,\vectsf{B}\,\text{sgn}(\vectsf{B}^\top\Tvect{z})
   -\Tvect{z}
   ,\label{eq::dis_solver_c_b}\\
 \dot{\tilde{\vect{x}}}&\!=\!-\vectsf{\Pi}\,\vect{\mathcal{H}}^{-1}\big(\Tvect{z}+\frac{\text{e}^{-t}}{N}\vect{1}_N\!\otimes\! \sum_{i=1}^N\vect{z}^i(0) \!+\!\,\vectsf{B}\vectsf{B}^\top{\Tvect{x}}\big)
 .
\end{align}
\end{subequations}
Here, we used~\eqref{eq::z_sum}. Moreover, given~\eqref{eq::change_of_var}  we obtain $\vect{1}_N\otimes \sum_{i=1}^N\tilde{\vect{z}}^i(t)=\vect{0}$  and $\vect{1}_N\otimes \sum_{i=1}^N\tilde{\vect{x}}^i(t)=\vect{0}$, which holds for any $t\in\real_{\geq0}$.

Next, we show that $\vectsf{B}^\top\tilde{\vect{z}}(t)$ goes to zero in finite time.
Defining $\hat{\vect{z}}=\vectsf{B}^\top\tilde{\vect{z}}(t)$ and $\hat{\vect{x}}=\vectsf{B}^\top\tilde{\vect{x}}(t)$, from~\eqref{eq::new_var} we get
\begin{subequations}\label{eq::new_var_hat}
\begin{align}
 \dot{\hat{\vect{z}}}&\!=\!- \vectsf{B}^\top\vectsf{B}\text{sgn}(\hat{\vect{z}})\!-\!\hat{\vect{z}}
 ,\label{eq::new_var_hat_a}\\
  \dot{\hat{\vect{x}}}&=\!-\vectsf{B}^\top\vect{\mathcal{H}}^{-1}(\Tvect{z}+\!\frac{\text{e}^{-t}}{N}\vect{1}_N\!\otimes\! \sum_{i=1}^N\vect{z}^i(0))-\!\,\vectsf{B}^\top\vect{\mathcal{H}}^{-1}\vectsf{B}\hat{\vect{x}}
  .\label{eq::new_var_hat_b} 
\end{align}
\end{subequations}
To analyze the stability of~\eqref{eq::new_var_hat_a}  consider  
\begin{align}\label{eq::Lya_func}
&V=\frac{1}{2}\hat{\vect{z}}^\top(\vectsf{B}^\top\vectsf{B})^{+}\hat{\vect{z}}\leq \bar{\lambda}\|\Hvect{z}\|_2^2,
\end{align}
where $\bar{\lambda}$ is the maximum eigenvalue of $\frac{1}{2}(\vectsf{B}^\top\vectsf{B})^{+}$. Note that $(\vectsf{B}^\top\vectsf{B})^{+}\geq0$. The Lie derivative of $V$ along~\eqref{eq::new_var_hat_a} is equal to
 \begin{align*}
 \dot{V}=&- \,\hat{\vect{z}}^\top\text{sgn}(\hat{\vect{z}})
 -\hat{\vect{z}}^\top\hat{\vect{z}}=- \|\Hvect{z}\|_1-\hat{\vect{z}}^\top\hat{\vect{z}}\leq- \|\Hvect{z}\|_1.
 \end{align*}
Since $\|\Hvect{z}\|_2\leq \|\Hvect{z}\|_1$, then we have 
$\dot{V}\leq-\|\hat{\vect{z}}\|_2\leq0$.
As such, from~\eqref{eq::Lya_func}, we obtain $\dot{V}\leq \frac{-1}{\sqrt{\bar{\lambda}}}\sqrt{V}$. Then, invoking the comparison Lemma~\cite[Lemma~3.4]{HKK:02}, we have the $\sqrt{V}\leq \sqrt{V(0)}-\frac{1}{2\sqrt{\bar{\lambda}}}t$.
Consequently, starting from any $V(0)$, $V(t)$ becomes zero at a finite time. Then, since $V=\frac{1}{2}\hat{\vect{z}}^\top(\vectsf{B}^\top\vectsf{B})^{+}\hat{\vect{z}}=\vect{z}^\top\vectsf{B}(\vectsf{B}^\top\vectsf{B})^{+}\vectsf{B}^\top{\vect{z}}=\vect{z}^\top\vectsf{\Pi}\vect{z}=(\vectsf{\Pi}\vect{z})^\top(\vectsf{\Pi}\vect{z})=\Tvect{z}^\top\Tvect{z}$,
and~\eqref{eq::z_g_sum}, we can conclude that every $\vect{z}^i$, $i\until{N}$, converges to $\frac{1}{N}\sum_{j=1}^N\vect{g}^j(\vect{x}^j)$ in finite time, and also $\Tvect{z}$ is bounded and converges to zero in finite time ($\|\Tvect{z}(t)\|_2\leq \|\Tvect{z}(0)\|_2-\frac{1}{2\sqrt{\bar{\lambda}}}t$
). 

Next, we show that $\vectsf{B}\hat{\vect{x}}$ converges to $0$ as $t\to\infty$. To this aim, we consider the radially unbounded Lyapunov function $W=\frac{1}{2}\hat{\vect{x}}^\top\hat{\vect{x}}$.
The Lie derivative of this function along the trajectories of~\eqref{eq::new_var_hat_b}~is
 \begin{align*}\label{eq::V_ddot}
 \dot{W}\!=\!&-\!\,\hat{\vect{x}}^\top\vectsf{B}^\top\!\vect{\mathcal{H}}^{-1}\vectsf{B}\hat{\vect{x}}
 \!-\!\hat{\vect{x}}^\top\vectsf{B}^\top\vect{\mathcal{H}}^{-1}\!(\tilde{\vect{z}}\!+\!\frac{\text{e}^{-t}}{N}\vect{1}_N\!\otimes\!\! \sum_{i=1}^N\!\!{\vect{z}}^i(0)).
 \end{align*}
Since $\Tvect{z}(t)$ vanishes in finite time, there exists a $t_1\in\real_{>0}$ such that for any $t\geq t_1$, we obtain
\begin{align*}
 &\dot{W}\!=-\,\|\sqrt{\vect{\mathcal{H}}^{-1}}\vectsf{B}\hat{\vect{x}}\|^2
 -\hat{\vect{x}}^\top\vectsf{B}^\top\vect{\mathcal{H}}^{-1}(\frac{\text{e}^{-t}}{N}\vect{1}_N\!\otimes\!\! \sum_{i=1}^N\!\!{\vect{z}}^i(0))\leq\\&\!\!-\!\|\!\sqrt{\!\vect{\mathcal{H}}^{-1}}\vectsf{B}\hat{\vect{x}}\|\big(\|\!\sqrt{\!\vect{\mathcal{H}}^{-1}}\vectsf{B}\hat{\vect{x}}\|\!-\!\frac{\text{e}^{-t}}{N}\|\!\sqrt{\!\vect{\mathcal{H}}^{-1}}(\vect{1}_N\!\otimes\!\! \sum_{i=1}^N\!\!{\vect{z}}^i(0))\|\big)
\end{align*}
Note that at each time $t\geq t_1$, $\dot{W}\leq0$ if and only if $\|\!\sqrt{\!\vect{\mathcal{H}}^{-1}}\vectsf{B}\hat{\vect{x}}\|\geq\frac{\text{e}^{-t}}{N}\|\sqrt{\!\vect{\mathcal{H}}^{-1}}(\vect{1}_N\!\otimes\!\! \sum_{i=1}^N\!\!{\vect{z}}^i(0))\|$. 
Then, at any $t_2\geq t_1$ if the trajectories of~\eqref{eq::new_var_hat} satisfy  $\Hvect{x}(t_2)\in\mathcal{S}=\{\hat{\vect{x}}\in\real^{Nd}\,|\,\,\|\sqrt{\vect{\mathcal{H}}^{-1}}\vectsf{B}\hat{\vect{x}}\|\leq\frac{\text{e}^{-t_2}}{\sqrt{\underline{m}}N}\|(\vect{1}_N\!\otimes\!\! \sum_{i=1}^N\!\!{\vect{z}}^i(0))\|\}$, then $\Hvect{x}(t)\in\mathcal{S}$ for all $t\geq t_2$, otherwise, $\dot{W}(t_2)< 0$. Here, we used the fact that $\|\sqrt{\vect{\mathcal{H}}^{-1}(\vect{x}(t))}\|\leq \frac{1}{\sqrt{\underline{m}}}$, which ensures that  $\frac{\text{e}^{-t}}{N}\|\sqrt{\vect{\mathcal{H}}^{-1}(\vect{x}(t))}(\vect{1}_N\!\otimes\!\! \sum_{i=1}^N\!\!{\vect{z}}^i(0))\|\leq\frac{\text{e}^{-t}}{\sqrt{\underline{m}}\,N}\|\vect{1}_N\!\otimes\!\! \sum_{i=1}^N\!\!{\vect{z}}^i(0)\|<\frac{\text{e}^{-t_2}}{\sqrt{\underline{m}}\,N} \|\vect{1}_N\!\otimes\!\! \sum_{i=1}^N\!\!{\vect{z}}^i(0)\|$ for any $t> t_2$. Therefore, as $t\to\infty$, we have the guarantees that along the trajectories of the system, $\|\sqrt{\vect{\mathcal{H}}(\vect{x}(t))^{-1}}\vectsf{B}\hat{\vect{x}}(t)\|$ goes to zero, which means that $\vectsf{B}\hat{\vect{x}}(t)$ goes to zero. Consequently, because of $\Hvect{x}=\vectsf{B}^\top\Tvect{x}$ and  $\vectsf{B}\vectsf{B}^\top=\vect{L}\otimes\vect{I}_d$ we can conclude that $(\vect{L}\otimes\vect{I}_d){\Tvect{x}}(t)$ goes to zero as $t\to\infty$. Given that the graph is connected, then, as $t\to\infty$, $\Hvect{x}$ goes to $\vect{1}_N\otimes \vect{\theta}$, $\vect{\theta}\in\real^d$, which given $\vect{1}_N\otimes \sum_{i=1}^N\tilde{\vect{x}}^i(t)=\vect{0}$, it means that $\Tvect{x}(t)$ goes to zero as $t\to\infty$. As a result, it follows from~\eqref{eq::change_of_var_b} that $\vect{x}^i$ converges to $\frac{1}{N}\sum_{j=1}^N \vect{x}^j(t)$ as $t\to\infty$. 
\end{pf}

\begin{rem}(Remark on the proof of Lemma~\ref{thm::convergence})\label{rem}
{\rm The solution of~\eqref{eq::dis_solver} is  in the sense of Filippov~\cite{AFF:88} since the solution is piecewise differentiable. However, the Filippov approach provides multi-valued functions for the solution  of~\eqref{eq::dis_solver} over discontinuity points; our stability analysis is valid since the Lyapunov function is smooth and decreasing over every Filippov solution of~\eqref{eq::dis_solver}.\boxend
}
\end{rem}

 Lemma~\ref{thm::convergence} showed that the trajectories of distributed HISO algorithm~\eqref{eq::dis_solver} converge to agreement space. The next theorem shows that this property indeed results in $\vect{x}^i$, $i\until{N}$ converging to $\vect{x}^\star$, the unique solution of the optimization problem~\eqref{eq::glob_opt_prob}, as $t\to\infty$.

\begin{thm}(Convergence of  algorithm~\eqref{eq::dis_solver})\label{thm::opt_solu}
{\rm Suppose the graph $\GG$ is connected and let Assumptions~\ref{asm:convexity} hold. Then, starting algorithm~\eqref{eq::dis_solver} over $\GG$ from any $\vect{x}^i(0),\vect{v}^i(0)\in\real^d$, $\sum_{i=1}^N\vect{v}^i(0)=\vect{0}$, every $\vect{x}^i(t)$, $i\until{N}$, converges  to $\vect{x}^\star$, the unique minimizer of~\eqref{eq::glob_opt_prob} and every $\vect{z}^i$ converges to zero as $t\rightarrow\infty$. 
}
\end{thm}
\begin{pf}
From~\eqref{eq::change_of_var_a}, $\vect{z}^i=\Tvect{z}^i+\frac{1}{N}\sum_{j=1}^N \vect{z}^i$. 
Recall from the proof of Lemma~\ref{thm::convergence} that under the stated initial condition, $\Tvect{z}^i$, $i\until{N}$, converges to zero in finite time. Then, convergence of $\vect{z}^i$, $i\!\until{N}$, to zero follows from~\eqref{eq::z_sum}. Next, notice  that~\eqref{eq::z_g_sum} and~\eqref{eq::z_sum} indicate that $\sum\nolimits_{j=1}^N\!\vect{g}^j(\vect{x}^j)$ goes to zero as $t\to\infty$. Then, because Lemma~\ref{thm::convergence} guarantees that $\vect{x}^i$, $i\!\until{N}$ converges to $\frac{1}{N}\sum_{j=1}^N\!\vect{x}^j$ as $t\!\to\!\infty$, we can conclude that $\frac{1}{N}\sum_{j=1}^N\!\vect{x}^j$, and subsequently every $\vect{x}^i$  converges to $\vect{x}^\star$ as~$t\!\to\!\infty$.
\end{pf}


\section{Numerical Example}\label{sec::Num_ex}
\begin{figure}[t]
  \centering
     \includegraphics[scale=0.68]{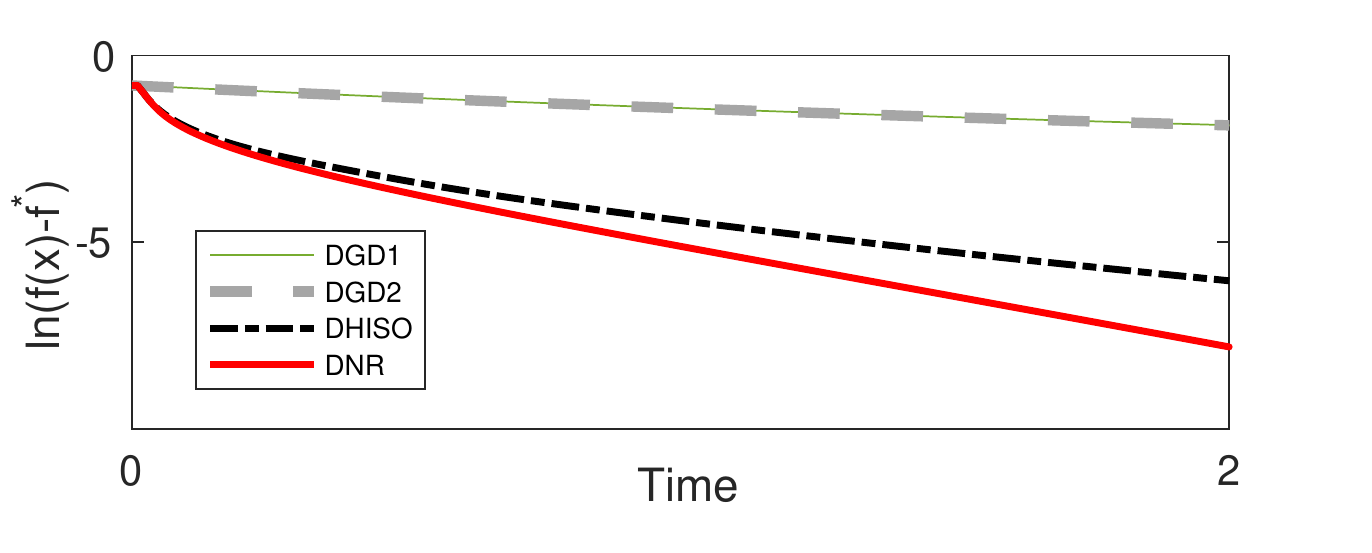}\vspace{-0.2in}
      \caption{ Convergence of DGD1, DGD2, DHISO and DNR of~\cite{DV-FZ-AC-GP-LS:15} algorithms in logaritmic scale.}\label{fig:example-1}
\end{figure}

We consider a distributed binary
classification problem using logistic regression over a connected graph of  Fig.~\ref{fig:network}. Each agent $i\until{N}$
has access to $m^i$ training samples
 $(\vect{c}_{ij} , y_{ij} )\in\real^{p}\times\{-1, +1\}$,
where $\vect{c}_{ij}$ contains p features of the $j^\text{th}$ training data at agent $i$,
and $y_{ij}$ is the corresponding binary label. The agents minimize 
$f=\sum_{i=1}^Nf^i(\vect{w},b)$ cooperatively, where $\vect{w}\in\real^p$, $b\in\real$, and each $f^i$
is given by  $f^i(\vect{w},b)=\sum\nolimits_{j=1}^{m^i} \ln(1+\ee^{-(\vect{w}^\top\vect{c}_{ij}+b)y_{ij}})+\frac{\lambda}{2N}\|\vect{w}\|^2.$
We generated the feature vectors $\vect{c}_{ij}$s randomly from two
distinct Gaussian distributions corresponding to two different labels, $+1$ and $-1$. Here, $p= 5$, $m = 10$, and $\lambda= 2$.
Figure~\ref{fig:example-1} shows the trajectories of the cost function   when the problem is solved via: distributed gradient descent algorithm of~\cite{SSK-JC-SM:15-auto} (DGD1), distributed gradient descent algorithm obtained from~\eqref{eq::dis_solver} when  $\vect{H}^i(\vect{x})$ are replaced by $\vect{I}_d$ (DGD2), our proposed distributed HISO algorithm (DHISO) and distributed NR algorithm (DNR) proposed in~\cite{DV-FZ-AC-GP-LS:15}.  As Fig.~\ref{fig:example-1} shows, DHISO and DNR algorithms both converge faster than  the gradient descent algorithms. Moreover, DHISO algorithm demonstrates a comparable response to that of the DNR but without requiring the neighboring agents to exchange their local Hessians with each other that the DNR algorithm of~\cite{DV-FZ-AC-GP-LS:15}~ requires.

\section{Conclusion}\label{sec::conclu}
We studied a novel second-order continuous-time distributed fast converging solution for an unconstrained optimization problem. Our approach guarantees  convergence to the  minimizer 
while keeping the communication cost efficient, in order of $O(Nd)$ as opposed to $O(Nd^2)$ for the existing results in the literature. Future work includes obtaining a discrete-time implementation of our algorithm with formal convergence guarantees.




\bibliographystyle{ieeetr}%
\bibliography{bib/alias,bib/Reference} 

\end{document}